\documentclass[12pt]{article}%
\usepackage{amsmath}
\usepackage{amsfonts}
\usepackage{amssymb}
\usepackage{graphicx}%
\begin{document}

\title{COMPUTING THE $\sin_{p}$ FUNCTION\\VIA THE INVERSE POWER METHOD}
\author{{\small \textbf{Rodney Josu\'{e} BIEZUNER, Grey ERCOLE, Eder Marinho MARTINS}%
}\thanks{ \textit{E-mail addresses:} rodney@mat.ufmg.br (R. J. Biezuner),
grey@mat.ufmg.br (G. Ercole), eder@iceb.ufop.br (E. Martins).}\\{\small \textit{Departamento de Matem\'{a}tica - ICEx, Universidade Federal de
Minas Gerais,}}\\{\small \textit{Av. Ant\^{o}nio Carlos 6627, Caixa Postal 702, 30161-970, Belo
Horizonte, MG, Brazil }} }
\maketitle

\noindent{\small ABSTRACT. In this paper, we discuss a new iterative method
for computing $\sin_{p}$. This function was introduced by Lindqvist in
connection with the unidimensional nonlinear Dirichlet eigenvalue problem for
the $p$-Laplacian. The iterative technique was inspired by the inverse power
method in finite dimensional linear algebra and is competitive with other
methods available in the literature.}


\section{Introduction}

In this paper we present a new method to compute the function $\sin_{p}$,
inspired by recent work done by the authors in \cite{BEM}, where an iterative
algorithm based on the inverse power method of linear algebra was introduced
for the computation of the first eigenvalue and first eigenfunction of the
Dirichlet problem for the $p$-Laplacian in arbitrary domains in $%
\mathbb{R}
^{N}$. The functions $\sin_{p}$, $1<p<\infty$, can be thought of as
generalizations of the familiar trigonometric functions. They arise in the
unidimensional Dirichlet eigenvalue problem for the $p$-Laplacian and were
introduced in this capacity in \cite{Lindqvist}, where a power series formula
for computing them was also formally given.

In \cite{BR1}\textbf{ }$\sin_{p}$\textbf{ }functions were utilized to
introduce a generalization of the Pr\"{u}fer transformation and thus
represent, in two phase-plane coordinates,\textbf{ }Sturm-Liouville-type
problems involving the $N$-dimensional radially symmetric $p$-Laplacian
$L_{p}u:=x^{1-N}\left(  x^{N-1}\left\vert u^{\prime}\right\vert ^{p-2}%
u^{\prime}\right)  ^{\prime}$, $0\leqslant a<x<b<\infty$.\textbf{ }This
approach was numerically implemented in \cite{BR2} for an eigenvalue problem
involving $L_{p}$ with separated homogeneous boundary conditions. In that
paper an interpolation table for $\sin_{p}$ was obtained by numerically
solving an ODE. Also in that paper the authors raised the question of finding
a fast and accurate algorithm for computing $\sin_{p}$.

Our method depends on the convergence of a sequence of functions whose
definition, as in \cite{BEM}, is motivated by an extension of the inverse
power method of linear algebra for obtaining the first eigenvalue and first
eigenfunction of finite dimensional linear operators. These functions are
recursively defined and can be given in integral form, so that they can be
obtained by numerical integration.

More specifically, recall that it suffices to obtain $\sin_{p}$ in the
interval $I_{p}=\left[  0,\pi_{p}/2\right]  $, since it is extended to the
interval $\left[  \pi_{p}/2,\pi_{p}\right]  $ symmetrically with respect to
$\pi_{p}/2$ and afterward to the whole real line $%
\mathbb{R}
$ as an odd, $2\pi_{p}$-periodic function (the definition of $\sin_{p}$ as
well as the precise value of $\pi_{p}$ are recalled in Section 2). We define
the following sequence of (positive) functions $\left\{  \phi_{n}\right\}
\subset C^{1}\left(  I_{p}\right)  $. Set $\phi_{0}\equiv1$ and
\[
\left\{
\begin{array}
[c]{ll}%
\left(  \phi_{n+1}^{\prime}\left\vert \phi_{n+1}^{\prime}\right\vert
^{p-2}\right)  ^{\prime}=-\phi_{n}\left\vert \phi_{n}\right\vert ^{p-2} &
\text{ \ \ if }x\in I_{p},\\
\phi_{n+1}\left(  0\right)  =\phi_{n+1}^{\prime}\left(  \pi_{p}/2\right)
=0. &
\end{array}
\right.
\]
We prove that the scaled sequence $\left\{  \sqrt[p]{p-1}\phi_{n}/\left\Vert
\phi_{n}\right\Vert _{\infty}\right\}  $ converges uniformly to $\sin_{p}$ in
$I_{p}$. The functions $\phi_{n}$ can be written in integral form as%
\[
\phi_{n+1}\left(  x\right)  =\int_{0}^{x}\left(  \int_{\theta}^{\pi_{p}/2}%
\phi_{n}\left(  s\right)  ^{p-1}ds\right)  ^{\frac{1}{p-1}}d\theta\text{,
\ \ }x\in I_{p},
\]
and, therefore, are readily computed using standard efficient numerical
methods for definite integrals.

This paper is organized as follows. In Section 2, we recall the definition and
some basic properties of $\sin_{p}$ which will be used in the sequel. In
Section 3, we show how to recursively construct a sequence of functions which
converge uniformly to $\sin_{p}$. Finally, in Section 4 we compare the
performance of our method with those of \cite{Lindqvist} and \cite{BR2}.

\section{The function $\sin_{p}$}

For the sake of completeness we recall in this section the definition and some
properties of the function $\sin_{p}$. The unidimensional Dirichlet eigenvalue
problem for the $p$-Laplacian, $p>1$, is%
\begin{equation}
\left\{
\begin{array}
[c]{ll}%
\psi_{p}\left(  u^{\prime}\right)  ^{\prime}=-\lambda\psi_{p}\left(  u\right)
& \quad\text{if }a<x<b,\\
u\left(  a\right)  =u\left(  b\right)  =0, &
\end{array}
\right.  \label{01}%
\end{equation}
where $\psi_{p}\left(  t\right)  =t\left\vert t\right\vert ^{p-2}$.

It is easy to verify that if $\lambda_{1}$ is the first eigenvalue of
\begin{equation}
\left\{
\begin{array}
[c]{ll}%
\psi_{p}\left(  v^{\prime}\right)  ^{\prime}=-\lambda\psi_{p}\left(  v\right)
& \quad\text{if }a<x<m:=\dfrac{a+b}{2},\\
v\left(  a\right)  =v^{\prime}\left(  m\right)  =0, &
\end{array}
\right.  \label{02}%
\end{equation}
and $v_{1}$ is the corresponding positive eigenfunction, then $\lambda_{1}$ is
also the first eigenvalue for (\ref{01}) with%
\[
u_{1}\left(  x\right)  =\left\{
\begin{array}
[c]{ll}%
v_{1}\left(  x\right)  & \quad\text{if }a\leqslant x\leqslant m,\\
v_{1}\left(  a+b-x\right)  & \quad\text{if }m\leqslant x\leqslant b,
\end{array}
\right.
\]
being the corresponding positive eigenfunction. Moreover, this function is
stricly increasing on $[a,m)$, strictly decreasing on $(m,b]$ and has only one
maximum point which is reached at $x=m$. Thus, $\left\Vert u_{1}\right\Vert
_{\infty}=u_{1}\left(  m\right)  $.

An expression for $\lambda_{1}$ is well known and can be obtained by
integration (see \cite{Otani}) as follows. First multiply (\ref{01}) by
$u_{1}^{\prime}$ and integrate the resulting equation by parts on $\left[
a,x\right]  $ to obtain%
\begin{equation}
\left.  \psi_{p}\left(  u_{1}^{\prime}\right)  u_{1}^{\prime}\right\vert
_{a}^{x}-\int_{a}^{x}\psi_{p}\left(  u_{1}^{\prime}\right)  u_{1}%
^{\prime\prime}dx=-\lambda_{1}\int_{a}^{x}\psi_{p}\left(  u_{1}\right)
u_{1}^{\prime}dx. \label{03}%
\end{equation}
We have%
\begin{equation}
\left.  \psi_{p}\left(  u_{1}^{\prime}\right)  u_{1}^{\prime}\right\vert
_{a}^{x}=\left\vert u_{1}^{\prime}\left(  x\right)  \right\vert ^{p}%
-\left\vert u_{1}^{\prime}\left(  a\right)  \right\vert ^{p} \label{04}%
\end{equation}%
\begin{equation}
\int_{a}^{x}\psi_{p}\left(  u_{1}\right)  u_{1}^{\prime}dx=\int_{u_{1}%
(a)}^{u_{1}(x)}\psi_{p}\left(  s\right)  \,ds=\frac{\left\vert u\left(
x\right)  \right\vert ^{p}}{p}-\frac{\left\vert u\left(  a\right)  \right\vert
^{p}}{p}, \label{05}%
\end{equation}%
\begin{equation}
\int_{a}^{x}\psi_{p}\left(  u_{1}^{\prime}\right)  u_{1}^{\prime\prime}%
dx=\int_{u_{1}^{\prime}\left(  a\right)  }^{u_{1}^{\prime}\left(  x\right)
}\psi_{p}\left(  s\right)  \,ds=\frac{\left\vert u^{\prime}\left(  x\right)
\right\vert ^{p}}{p}-\frac{\left\vert u^{\prime}\left(  a\right)  \right\vert
^{p}}{p}. \label{06}%
\end{equation}
Substituting (\ref{04}), (\ref{05}) and (\ref{06}) in (\ref{03}) we obtain
\[
\left(  1-\frac{1}{p}\right)  \left[  \left\vert u_{1}^{\prime}\left(
x\right)  \right\vert ^{p}-\left\vert u_{1}^{\prime}\left(  a\right)
\right\vert ^{p}\right]  =-\lambda_{1}\left[  \frac{\left\vert u_{1}\left(
x\right)  \right\vert ^{p}}{p}-\frac{\left\vert u_{1}\left(  a\right)
\right\vert ^{p}}{p}\right]  ,
\]
whence
\[
\left[  \left(  1-\frac{1}{p}\right)  \left\vert u_{1}^{\prime}\right\vert
^{p}+\lambda_{1}\frac{\left\vert u_{1}\right\vert ^{p}}{p}\right]  _{a}%
^{x}=0.
\]
This means that
\[
\frac{p-1}{p}\left\vert u_{1}^{\prime}\right\vert ^{p}+\frac{\lambda_{1}}%
{p}\left\vert u_{1}\right\vert ^{p}\equiv C,
\]
where $C$ is a constant and $p^{\prime}=p/\left(  p-1\right)  $ is the
conjugate of $p$. The value of $C$ can be found computing the value of this
expression at the maximum point $m$; choosing $u_{1}$ such that $u_{1}\left(
m\right)  =1$ we find%
\[
C=\frac{p-1}{p}\left\vert u_{1}^{\prime}\left(  m\right)  \right\vert
^{p}+\frac{\lambda_{1}}{p}\left\vert u_{1}\left(  m\right)  \right\vert
^{p}=\frac{\lambda_{1}}{p}.
\]
Therefore,%
\begin{equation}
\left(  p-1\right)  \left\vert u_{1}^{\prime}\left(  x\right)  \right\vert
^{p}+\lambda_{1}\left\vert u_{1}\left(  x\right)  \right\vert ^{p}=\lambda_{1}
\label{07}%
\end{equation}
for all $x\in\left[  a,b\right]  $.

On the interval $\left[  a,m\right]  $ we have $u^{\prime}\geqslant0$, hence
we can write%
\begin{equation}
\frac{u_{1}^{\prime}\left(  x\right)  }{\sqrt[p]{\left(  1-\left\vert
u_{1}\left(  x\right)  \right\vert ^{p}\right)  }}=\sqrt[p]{\frac{\lambda_{1}%
}{p-1}} \label{08}%
\end{equation}
for all $x\in\left[  a,m\right]  $. Integrating this equation on $\left(
a,m\right)  $ leads to%
\[
\frac{b-a}{2}\sqrt[p]{\frac{\lambda_{1}}{p-1}}=\int_{u_{1}(a)}^{u_{1}(m)}%
\frac{ds}{\sqrt[p]{1-s^{p}}}=\int_{0}^{1}\frac{ds}{\sqrt[p]{1-s^{p}}},
\]
which gives the expression%
\begin{equation}
\lambda_{1}=(p-1)\left(  \frac{2}{b-a}\int_{0}^{1}\frac{ds}{\sqrt[p]{1-s^{p}}%
}\right)  ^{p}=\left(  \frac{\pi_{p}}{b-a}\right)  ^{p}, \label{09}%
\end{equation}
where we set%
\begin{equation}
\pi_{p}:=2\sqrt[p]{p-1}\int_{0}^{1}\frac{ds}{\sqrt[p]{1-s^{p}}}. \label{ppip}%
\end{equation}
Making the change of variable $s=\sqrt[p]{t}$ in the last integral and using
the classical Beta function $B$ we obtain
\[
\int_{0}^{1}\frac{ds}{\sqrt[p]{1-s^{p}}}=\frac{1}{p}\int_{0}^{1}t^{\frac{1}%
{p}-1}(1-t)^{-\frac{1}{p}}dt=\frac{1}{p}B\left(  1-\frac{1}{p},\frac{1}%
{p}\right)  =\frac{\pi/p}{\sin(\pi/p)}%
\]
(Here one use the properties $B(x,y)B(x+y,1-y)=x/x\sin\left(  \pi y\right)  $
and $B(1,z)=1/z$ with $x=1-1/p$ and $y=z=1/p$).

Therefore,%
\begin{equation}
\pi_{p}=\frac{2\sqrt[p]{p-1}\left(  \pi/p\right)  }{\sin(\pi/p)} \label{10}%
\end{equation}
and
\[
\lambda_{1}=\left(  \frac{2\sqrt[p]{p-1}\left(  \pi/p\right)  }{(b-a)\sin
\left(  \pi/p\right)  }\right)  ^{p}.
\]
When $a=0$ and $b=\pi_{p}$ we denote the function $\sqrt[p]{p-1}u_{1}$ by
$\sin_{p}.$ Thus, $\sin_{p}\left(  0\right)  =0=\sin_{p}^{\prime}\left(
\pi_{p}/2\right)  $, $\lambda_{1}=1$ and from (\ref{07}):%
\[
\left\vert \sin_{p}^{\prime}\right\vert ^{p}+\frac{\left\vert \sin
_{p}\right\vert ^{p}}{p-1}=1.
\]
It is clear from this equation that $\sin_{p}^{\prime}\left(  0\right)  =1.$

We remark that $u=\sin_{p}$ is also the unique solution of the initial value
problem
\[
\left\vert u^{\prime}\right\vert ^{p}+\frac{\left\vert u\right\vert ^{p}}%
{p-1}=1,\quad u\left(  0\right)  =0,
\]
which can be used to define this function.

Alternatively, we can define $\sin_{p}$ on the interval $\left[  0,\pi
_{p}/2\right]  $ as an inverse function. In fact, multiplying (\ref{08}) by
$\sqrt[p]{p-1}$ and using (\ref{09}) with $a=0$ and $b=\pi_{p}$ we obtain
\[
\int_{0}^{\sin_{p}\left(  x\right)  }\frac{ds}{\sqrt[p]{\left(  1-\frac{s^{p}%
}{p-1}\right)  }}=x,\quad\text{for }x\in\left[  0,\pi_{p}/2\right]  ,
\]
that is, $\sin_{p}=\zeta^{-1}$ where%
\[
\zeta\left(  z\right)  :=%
{\displaystyle\int_{0}^{z}}
\frac{ds}{\sqrt[p]{\left(  1-\frac{s^{p}}{p-1}\right)  }},\quad\text{for }%
z\in\left[  0,\sqrt[p]{p-1}\right]  .
\]

With this definition, we extend $\sin_{p}$ to the interval $\left[  \pi
_{p}/2,\pi_{p}\right]  $ symmetrically with respect to $\pi_{p}/2$ and
afterward to the whole real line $%
\mathbb{R}
$ as an odd, $2\pi_{p}$-periodic function. We list the basic properties of
$\sin_{p}$:

\begin{enumerate}
\item $\sin_{p}\left(  0\right)  =0=\sin_{p}\left(  \pi_{p}\right)  $,
$\sin_{p}\left(  \pi_{p}/2\right)  =\left\Vert \sin_{p}\right\Vert _{\infty
}=\sqrt[p]{p-1}$.

\item $\sin_{p}\left(  x\right)  $ is strictly increasing in $\left[
0,\pi_{p}/2\right]  $ and strictly decreasing in $\left[  \pi_{p}/2,\pi
_{p}\right]  .$

\item $\left\vert \sin_{p}^{\prime}\left(  x\right)  \right\vert
=\sqrt[p]{1-\dfrac{\left\vert \sin_{p}\right\vert ^{p}}{p-1}}.$
\end{enumerate}

\section{A sequence uniformly convergent to $\sin_{p}$}

Let $I_{p}=\left[  0,\pi_{p}/2\right]  $ and define the following sequence of
functions $\left\{  \phi_{n}\right\}  \subset C^{1}\left(  I_{p}\right)  $.
Set $\phi_{0}\equiv1$ and
\[
\left\{
\begin{array}
[c]{ll}%
\left(  \psi_{p}\left(  \phi_{n+1}^{\prime}\right)  \right)  ^{\prime}%
=-\psi_{p}\left(  \phi_{n}\right)  & \text{ \ \ if }x\in I_{p},\\
\phi_{n+1}\left(  0\right)  =\phi_{n+1}^{\prime}\left(  \pi_{p}/2\right)
=0. &
\end{array}
\right.
\]
In this section, we prove that the scaled sequence $\left\{  \sqrt[p]{p-1}%
\phi_{n}/\left\Vert \phi_{n}\right\Vert _{\infty}\right\}  $ converges
uniformly to $\sin_{p}$ in $I_{p}$. Before proceeding, we recall some basic
properties of the $\psi_{p}$ functions:

\begin{description}
\item[Proposition 3.1.] (Basic properties of $\psi_{p}$) \textit{The following
holds:}
\end{description}

\begin{enumerate}
\item $\psi_{p}$\textit{ is continuous, strictly increasing and odd, for each
}$p>1.$

\item $\psi_{p}\left(  ab\right)  =\psi_{p}\left(  a\right)  \psi_{p}\left(
b\right)  .$

\item $\psi_{p}\left(  \dfrac{a}{b}\right)  =\dfrac{\psi_{p}\left(  a\right)
}{\psi_{p}\left(  b\right)  }$

\item $\left(  \psi_{p}\right)  ^{-1}=\psi_{p^{\prime}}.$

\item $\int_{0}^{t}\psi_{p}\left(  s\right)  ds=\dfrac{\left\vert t\right\vert
^{p}}{p}.$
\end{enumerate}

By a straightforward calculation we can find the following recursive integral
expression for the $\phi_{n}$-functions:
\begin{equation}
\phi_{n+1}\left(  x\right)  =\int_{0}^{x}\psi_{p^{\prime}}\left(  \int
_{\theta}^{\pi_{p}/2}\psi_{p}\left(  \phi_{n}\left(  s\right)  \right)
ds\right)  d\theta. \label{11}%
\end{equation}
It is clear from (\ref{11}) that each $\phi_{n}\ $is positive, increasing on
$I_{p}$ and reaches its maximum value at $x=\pi_{p}/2$. One can obtain an
explicit expression for $\phi_{1}$, the second function in the sequence:
\begin{align*}
\phi_{1}\left(  x\right)   &  =\int_{0}^{x}\psi_{p^{\prime}}\left(
\int_{\theta}^{\pi_{p}/2}\psi_{p}\left(  1\right)  ds\right)  d\theta\\
&  =\int_{0}^{x}\psi_{p^{\prime}}\left(  \frac{\pi_{p}}{2}-\theta\right)
d\theta\\
&  =\int_{\pi_{p}/2-x}^{\pi_{p}/2}\psi_{p^{\prime}}\left(  y\right)  dy\\
&  =\frac{1}{p}\left[  \left(  \frac{\pi_{p}}{2}\right)  ^{p}-\left(
\frac{\pi_{p}}{2}-x\right)  ^{p}\right]  .
\end{align*}
Note that%
\[
\left\Vert \phi_{1}\right\Vert _{\infty}=\phi_{1}\left(  \frac{\pi_{p}}%
{2}\right)  =\frac{1}{p}\left(  \frac{\pi_{p}}{2}\right)  ^{p}=\frac{p-1}%
{p}\left(  \frac{\pi/p}{\sin\left(  \pi/p\right)  }\right)  ^{p}.
\]
The next $\phi_{n}$-functions however, are very difficult to obtain explicitly
by solving the integrals analytically. On the other hand, the integrals can
easily be solved numerically.

\begin{description}
\item[Proposition 3.2.] $\phi_{n+1}\leqslant\left\Vert \phi_{1}\right\Vert
_{\infty}\phi_{n}$\textit{ on }$I_{p}.$
\end{description}

\noindent\textbf{Proof.} For $n=1$ the result is trivially true since
$\phi_{0}\equiv1.$Assuming by induction that $\phi_{n}\leqslant\left\Vert
\phi_{1}\right\Vert _{\infty}\phi_{n-1}$, we have
\begin{align*}
\phi_{n+1}\left(  x\right)   &  =\int_{0}^{x}\psi_{p^{\prime}}\left(
\int_{\theta}^{\pi_{p}/2}\psi_{p}\left(  \phi_{n}\left(  s\right)  \right)
\,ds\right)  \,d\theta\\
&  \leqslant\int_{0}^{x}\psi_{p^{\prime}}\left(  \int_{\theta}^{\pi_{p}/2}%
\psi_{p}\left(  \left\Vert \phi_{1}\right\Vert _{\infty}\phi_{n-1}\left(
s\right)  \right)  \,ds\right)  \,d\theta\\
&  =\int_{0}^{x}\psi_{p^{\prime}}\left(  \psi_{p}\left(  \left\Vert \phi
_{1}\right\Vert _{\infty}\right)  \int_{\theta}^{\pi_{p}/2}\psi_{p}\left(
\phi_{n-1}\left(  s\right)  \right)  \,ds\right)  \,d\theta\\
&  =\left\Vert \phi_{1}\right\Vert _{\infty}\int_{0}^{x}\psi_{p^{\prime}%
}\left(  \int_{\theta}^{\pi_{p}/2}\psi_{p}\left(  \phi_{n-1}\left(  s\right)
\right)  \,ds\right)  \,d\theta\\
&  =\left\Vert \phi_{1}\right\Vert _{\infty}\phi_{n}\left(  x\right)  .
\end{align*}

\noindent$\blacksquare$

The following technical lemma, which will be used in the sequel, can be proved
via the Cauchy mean value theorem (see \cite{AVV}) and works as a
L'H\^{o}pital's rule in order to get monotonicity for a certain quotient function.

\begin{description}
\item[Lemma 3.3.] \textit{Let }$f,g:\left[  a,b\right]  \longrightarrow
R$\textit{ be continuous on }$\left[  a,b\right]  $\textit{ and differentiable
in }$\left(  a,b\right)  $. \textit{Suppose }$g^{\prime}(x)\neq0$\textit{ for
all }$x\in\left(  a,b\right)  $. \textit{If }$\dfrac{f^{\prime}}{g^{\prime}}%
$\textit{ is (strictly) increasing }[\textit{decreasing}],\textit{ then both
}$\dfrac{f(x)-f(a)}{g(x)-g(a)}$\textit{ and }$\dfrac{f(x)-f(b)}{g(x)-g(b)}%
$\textit{are (strictly) increasing }[\textit{decreasing}]\textit{.}

\item[Theorem 3.4.] \textit{For each }$n\geqslant1$\textit{ the function
}$\dfrac{\phi_{n}}{\phi_{n+1}}$\textit{ is strictly decreasing on }$I_{p}%
$\textit{ and}
\end{description}

\begin{enumerate}
\item[(i)] $\dfrac{1}{\left\Vert \phi_{1}\right\Vert _{\infty}}\leqslant
\inf\limits_{I_{p}}\dfrac{\phi_{n}}{\phi_{n+1}}=\dfrac{\phi_{n}\left(  \pi
_{p}/2\right)  }{\phi_{n+1}\left(  \pi_{p}/2\right)  }=\dfrac{\left\Vert
\phi_{n}\right\Vert _{\infty}}{\left\Vert \phi_{n+1}\right\Vert _{\infty}}.$

\item[(ii)] $\left\Vert \dfrac{\phi_{n}}{\phi_{n+1}}\right\Vert _{\infty}%
=\psi_{p^{\prime}}\left(  \dfrac{\int_{0}^{\pi_{p}/2}\psi_{p}\left(
\phi_{n-1}\left(  s\right)  \right)  ds}{\int_{0}^{\pi_{p}/2}\psi_{p}\left(
\phi_{n}\left(  s\right)  \right)  ds}\right)  $ \ \ for $n\geqslant1.$

\item[(iii)] $\left\Vert \dfrac{\phi_{n}}{\phi_{n+1}}\right\Vert _{\infty
}\leqslant\left\Vert \dfrac{\phi_{n-1}}{\phi_{n}}\right\Vert _{\infty
}\leqslant\cdots\leqslant\left\Vert \dfrac{\phi_{1}}{\phi_{2}}\right\Vert
_{\infty}<\infty.$
\end{enumerate}

\noindent\textbf{Proof.} Since $\phi_{1}$ is strictly increasing, it follows
that $1/\phi_{1}$ is strictly decreasing. Assume by induction that $\phi
_{n-1}/\phi_{n}$ is strictly decreasing. Since
\[
\frac{\phi_{n}\left(  x\right)  -\phi_{n}\left(  0\right)  }{\phi_{n+1}%
-\phi_{n+1}\left(  0\right)  }=\frac{\phi_{n}\left(  x\right)  }{\phi
_{n+1}\left(  x\right)  },
\]
in order to show that $\phi_{n}/\phi_{n+1}$ is strictly decreasing, it
suffices in light of the lemma to verify that $\phi_{n}^{\prime}/\phi
_{n+1}^{\prime}$ is strictly decreasing on $I_{p}$. But,
\[
\frac{\phi_{n}^{\prime}\left(  x\right)  }{\phi_{n+1}^{\prime}\left(
x\right)  }=\dfrac{\psi_{p^{\prime}}\left(
{\displaystyle\int_{x}^{\pi_{p}/2}}
\psi_{p}\left(  \phi_{n-1}\left(  s\right)  \right)  ds\right)  }%
{\psi_{p^{\prime}}\left(
{\displaystyle\int_{x}^{\pi_{p}/2}}
\psi_{p}\left(  \phi_{n}\left(  s\right)  \right)  ds\right)  }=\psi
_{p^{\prime}}\left(  \frac{%
{\displaystyle\int_{x}^{\pi_{p}/2}}
\psi_{p}\left(  \phi_{n-1}\left(  s\right)  \right)  ds}{%
{\displaystyle\int_{x}^{\pi_{p}/2}}
\psi_{p}\left(  \phi_{n}\left(  s\right)  \right)  ds}\right)  .
\]
Since $\psi_{p^{\prime}}$ is strictly increasing and the functions $\int
_{x}^{\pi_{p}/2}\psi_{p}\left(  \phi_{n-1}\left(  s\right)  \right)  ds$ and
$\int_{x}^{\pi_{p}/2}\psi_{p}\left(  \phi_{n}\left(  s\right)  \right)  ds$
are null at $x=\pi_{p}/2$, we can apply the lemma again to verify that the
quotient of these integral functions is a strictly decreasing function. We
have%
\[
\frac{\left(
{\displaystyle\int_{x}^{\pi_{p}/2}}
\psi_{p}\left(  \phi_{n-1}\left(  s\right)  \right)  ds\right)  ^{\prime}%
}{\left(
{\displaystyle\int_{x}^{\pi_{p}/2}}
\psi_{p}\left(  \phi_{n}\left(  s\right)  \right)  ds\right)  ^{\prime}}%
=\frac{\psi_{p}\left(  \phi_{n-1}\left(  s\right)  \right)  }{\psi_{p}\left(
\phi_{n}\left(  s\right)  \right)  }=\psi_{p}\left(  \frac{\phi_{n-1}}%
{\phi_{n}}\right)  ,
\]
which is strictly decreasing by the induction hypothesis.

The inequality in (i) follows from Proposition 3.2. Before verifying (ii) we
remark that $\left\Vert 1/\phi_{1}\right\Vert _{\infty}=\infty$ since
$\phi_{1}\left(  0\right)  =0.$ In order to prove (ii) we first observe that
the monotonicity of $\phi_{n}/\phi_{n+1}$ implies that
\[
\left\Vert \dfrac{\phi_{n}}{\phi_{n+1}}\right\Vert _{\infty}=\lim
_{x\rightarrow0^{+}}\frac{\phi_{n}\left(  x\right)  }{\phi_{n+1}\left(
x\right)  }.
\]
L'H\^{o}pital's rule then yields%
\[
\lim_{x\rightarrow0^{+}}\frac{\phi_{n}\left(  x\right)  }{\phi_{n+1}\left(
x\right)  }=\lim_{x\rightarrow0^{+}}\frac{\phi_{n}^{\prime}\left(  x\right)
}{\phi_{n+1}^{\prime}\left(  x\right)  }=\psi_{p^{\prime}}\left(  \frac
{\int_{0}^{\pi_{p}/2}\psi_{p}\left(  \phi_{n-1}\left(  s\right)  \right)
ds}{\int_{0}^{\pi_{p}/2}\psi_{p}\left(  \phi_{n}\left(  s\right)  \right)
ds}\right)  <\infty.
\]
The proof of (iii) is a consequence of the following estimates, valid for
$n\geqslant2$:%
\begin{align*}
\left\Vert \frac{\phi_{n}}{\phi_{n+1}}\right\Vert _{\infty}  &  =\psi
_{p^{\prime}}\left(  \frac{%
{\displaystyle\int_{0}^{\pi_{p}/2}}
\psi_{p}\left(  \phi_{n-1}\left(  s\right)  \right)  ds}{%
{\displaystyle\int_{0}^{\pi_{p}/2}}
\psi_{p}\left(  \phi_{n}\left(  s\right)  \right)  ds}\right) \\
&  \leqslant\psi_{p^{\prime}}\left(  \frac{%
{\displaystyle\int_{0}^{\pi_{p}/2}}
\psi_{p}\left(  \phi_{n}\left(  s\right)  \right)  \psi_{p}\left(  \dfrac
{\phi_{n-1}}{\phi_{n}}\left(  s\right)  \right)  ds}{%
{\displaystyle\int_{0}^{\pi_{p}/2}}
\psi_{p}\left(  \phi_{n}\left(  s\right)  \right)  ds}\right) \\
&  \leqslant\psi_{p^{\prime}}\left(  \frac{%
{\displaystyle\int_{0}^{\pi_{p}/2}}
\psi_{p}\left(  \phi_{n}\left(  s\right)  \right)  \psi_{p}\left(  \left\Vert
\dfrac{\phi_{n-1}}{\phi_{n}}\right\Vert _{\infty}\right)  ds}{%
{\displaystyle\int_{0}^{\pi_{p}/2}}
\psi_{p}\left(  \phi_{n}\left(  s\right)  \right)  ds}\right) \\
&  =\left\Vert \frac{\phi_{n-1}}{\phi_{n}}\right\Vert _{\infty}\psi
_{p^{\prime}}\left(  \frac{%
{\displaystyle\int_{0}^{\pi_{p}/2}}
\psi_{p}\left(  \phi_{n}\left(  s\right)  \right)  ds}{%
{\displaystyle\int_{0}^{\pi_{p}/2}}
\psi_{p}\left(  \phi_{n}\left(  s\right)  \right)  ds}\right) \\
&  =\left\Vert \frac{\phi_{n-1}}{\phi_{n}}\right\Vert _{\infty}.
\end{align*}

\noindent$\blacksquare$

\begin{description}
\item[Theorem 2.4.] \textit{Let }$u_{n}:=\dfrac{\phi_{n}}{\left\Vert \phi
_{n}\right\Vert _{\infty}}\in C^{1}\left(  I_{p}\right)  ,$\textit{ for
}$n\geqslant1.$\textit{ Then the sequence }$\left\{  u_{n}\left(  x\right)
\right\}  _{n\geqslant1}$\textit{ is decreasing for each }$x\in I_{p}$\textit{
and }%
\[
\sqrt[p]{p-1}u_{n}\rightarrow\sin_{p}\text{\textit{ \ \ uniformly in} }I_{p}.
\]

\end{description}

\noindent\textbf{Proof.} In $I_{p}$ we have%
\begin{align*}
\frac{u_{n}}{u_{n+1}}  &  =\dfrac{\phi_{n}}{\phi_{n+1}}\left(  \dfrac
{\left\Vert \phi_{n}\right\Vert _{\infty}}{\left\Vert \phi_{n+1}\right\Vert
_{\infty}}\right)  ^{-1}\\
&  \geqslant\left(  \inf\limits_{I_{p}}\dfrac{\phi_{n}}{\phi_{n+1}}\right)
\left(  \dfrac{\left\Vert \phi_{n}\right\Vert _{\infty}}{\left\Vert \phi
_{n+1}\right\Vert _{\infty}}\right)  ^{-1}\\
&  =\left(  \dfrac{\left\Vert \phi_{n}\right\Vert _{\infty}}{\left\Vert
\phi_{n+1}\right\Vert _{\infty}}\right)  \left(  \dfrac{\left\Vert \phi
_{n}\right\Vert _{\infty}}{\left\Vert \phi_{n+1}\right\Vert _{\infty}}\right)
^{-1}\\
&  =1,
\end{align*}
that is, $\left\{  u_{n}\left(  x\right)  \right\}  _{n\geqslant1}$ is
decreasing for each $x\in I_{p},$ and the whole sequence is bounded below by
$u_{1}$. Thus, there exists%
\[
u:=\lim u_{n}.
\]
We have $\left\Vert u_{n}\right\Vert _{\infty}=1$ for each $n$. Moreover,
since%
\[
\frac{\left\Vert \phi_{n}\right\Vert _{\infty}}{\left\Vert \phi_{n+1}%
\right\Vert _{\infty}}=\inf\limits_{I_{p}}\dfrac{\phi_{n}}{\phi_{n+1}%
}\leqslant\left\Vert \frac{\phi_{n}}{\phi_{n+1}}\right\Vert _{\infty}%
\leqslant\left\Vert \frac{\phi_{1}}{\phi_{2}}\right\Vert _{\infty}=:C,
\]
we also have, for every $x\in I_{p},$%
\begin{align*}
\left\vert u_{n}^{\prime}\left(  x\right)  \right\vert  &  =\frac
{1}{\left\Vert \phi_{n}\right\Vert _{\infty}}\psi_{p^{\prime}}\left(  \int
_{x}^{\pi_{p}/2}\psi_{p}\left(  \phi_{n-1}\left(  s\right)  \right)  ds\right)
\\
&  =\frac{\left\Vert \phi_{n-1}\right\Vert _{\infty}}{\left\Vert \phi
_{n}\right\Vert _{\infty}}\psi_{p^{\prime}}\left(  \int_{x}^{\pi_{p}/2}%
\psi_{p}\left(  \frac{\phi_{n-1}\left(  s\right)  }{\left\Vert \phi
_{n-1}\right\Vert _{\infty}}\right)  ds\right) \\
&  \leqslant C\psi_{p^{\prime}}\left(  \int_{0}^{\pi_{p}/2}\psi_{p}\left(
u_{n-1}\right)  ds\right) \\
&  \leqslant C\psi_{p^{\prime}}\left(  \int_{0}^{\pi_{p}/2}\psi_{p}\left(
1\right)  ds\right) \\
&  =\frac{C\pi_{p}}{2}.
\end{align*}
It follows from Arzela-Ascoli's theorem that $u_{n}\rightarrow u\in C\left(
I_{p}\right)  $, uniformly.

In order to conclude the proof, we need just to show that
\begin{equation}
u=\dfrac{\sin_{p}}{\sqrt[p]{p-1}}. \label{12}%
\end{equation}
From (\ref{11}) we can write the following expression:
\[
u_{n+1}\left(  x\right)  =\gamma_{n}\int_{0}^{x}\psi_{p^{\prime}}\left(
\int_{\theta}^{\pi_{p}/2}\psi_{p}\left(  u_{n}\left(  s\right)  \right)
ds\right)  d\theta,
\]
where
\[
\gamma_{n}:=\frac{\left\Vert \phi_{n}\right\Vert _{\infty}}{\left\Vert
\phi_{n+1}\right\Vert _{\infty}}.
\]
In view of the boundedness of $\left\{  \gamma_{n}\right\}  $, there exists
$\gamma:=\lim\gamma_{n_{k}}$ for some subsequence $\left\{  \gamma_{n_{k}%
}\right\}  .$ Thus, letting $k\rightarrow\infty$ in
\[
u_{n_{k}+1}\left(  x\right)  =\gamma_{n_{k}}\int_{0}^{x}\psi_{p^{\prime}%
}\left(  \int_{\theta}^{\pi_{p}/2}\psi_{p}\left(  u_{n_{k}}\left(  s\right)
\right)  ds\right)  d\theta,
\]
we get%
\[
u\left(  x\right)  =\gamma\int_{0}^{x}\psi_{p^{\prime}}\left(  \int_{\theta
}^{\pi_{p}/2}\psi_{p}\left(  u\left(  s\right)  \right)  ds\right)  d\theta\in
C^{1}\left(  I_{p}\right)  ,
\]
which means that $u$ is a positive solution to the following problem%
\[
\left\{
\begin{array}
[c]{ll}%
\psi_{p}\left(  u^{\prime}\right)  ^{\prime}=-\gamma\psi_{p}\left(  u\right)
& \text{ \ \ if }x\in I_{p},\\
u\left(  0\right)  =u^{\prime}\left(  \pi_{p}/2\right)  =0. &
\end{array}
\right.
\]
In view of the positivity of $u$, we can integrate the equation above
multiplied by $u^{\prime}$ and proceed as in the derivation of (\ref{09}) to
find $\gamma=1$. From this we conclude that in fact $\lim\gamma_{n}=1$ (the
whole sequence converges to the eingenvalue $1$) and that $u=\lim u_{n}$
satisfies the same boundary value problem that $\sin_{p}/\left\Vert \sin
_{p}\right\Vert $ does. Since both $u$ and $\sin_{p}$ are positive and
$\left\Vert u\right\Vert _{\infty}=\left\Vert \sin_{p}/\left\Vert \sin
_{p}\right\Vert \right\Vert _{\infty}=1$, we must have
\[
u=\dfrac{\sin_{p}}{\left\Vert \sin_{p}\right\Vert }%
\]
whence (\ref{12}) follows. $\blacksquare$

\section{Numerical Results}

Next we examine the computational time of each method. Computations were
performed on a WindowsXP/Pentium 4-2.8GHz platform, using the GCC compiler.
Although the method of computing $\sin_{p}$ by solving an ODE suggested in
\cite{BR2} (which we implemented by means of a standard Runge-Kutta fourth
power method) is by far the fastest, the computational times of the other two
methods are competitive, the inverse power method being on average more than
twice as fast as the power series method of \cite{Lindqvist} for values of $p$
greater than 2. Also, the average number of $8$ iterations that the inverse
power method uses to obtain the same (and sometimes better; see Table 2)
accuracy of the differential equation method of \cite{BR2} is quite
remarkable, specially taking into account that the functions $\phi_{n}$
converge to $0$ rather rapidly. We emphasize that the computational time of
the inverse power method is not the main subject of this presentation. The
method demands the computation of double integrals at each iteration for each
grid point. We opted for a classical, computationally easy to implement and
reasonably fast method to compute these integrals, namely, the Simpson
composite method. However, a greater effort spent in lessening the
computational time of the numerical integrations certainly would be reflected
in a substantial decrease in the time spent computing $\sin_{p}$ overall.
Nevertheless, by considering the accuracy and the comparison scale among the
three methods (on the range of miliseconds) we may say that the results
presented in this paper validate the inverse power method as an effective and
reasonably fast method for numerically obtaining $\sin_{p}$.

Below we present the average time spent in computing $\sin_{p}$ on the whole
interval $I_{p}$ divided in $101$ grid points by each method for six values of
$p$ (the average was taken out of five computer runs); the stop criterion in
each method was an error tolerance of $10^{-8}$ between successive iterations
and less than 500 terms in the power series.%

\[%
\begin{tabular}
[c]{|l|c|c|c|c|c|c|}\hline
$p$ & $1.1$ & $1.5$ & $2.0$ & $2.5$ & $3.0$ & $3.5$\\\hline
Inverse power method & $21.5$ & $32.1$ & $1.1$ & $37.7$ & $37.8$ &
$31.7$\\\hline
Differential equation method & $1.9$ & $1.8$ & $1.1$ & $1.5$ & $1.5$ &
$1.5$\\\hline
Power series & $92.9$ & $2.2$ & $2.0$ & $79.6$ & $79.3$ & $73.3$\\\hline
\end{tabular}
\ \
\]

\begin{center}
{\small Table 1: Average time (in miliseconds) for the computation of
$\sin_{p}$ on $I_{p}$ for each method.}\vspace{0.2cm}
\end{center}

Besides the trivial point $0$, the only point where the value of $\sin_{p}$ is
exactly known is $\pi_{p}/2$, with $\sin_{p}\left(  \pi_{p}/2\right)
=\sqrt[p]{p-1}$. In the next table we present the computed value for $\sin
_{p}\left(  \pi_{p}/2\right)  $ obtained using each method:%

\[%
\begin{tabular}
[c]{|l|c|c|c|c|c|c|}\hline
$p$ & $1.1$ & $1.5$ & $2.0$ & $2.5$ & $3.0$ & $3.5$\\\hline
$\sqrt[p]{p-1}$ & $0.123285$ & $0.629961$ & $1$ & $1.17608$ & $1.25992$ &
$1.29926$\\\hline
Inverse power method & $0.123285$ & $0.629961$ & $1$ & $1.17608$ & $1.25992$ &
$1.29926$\\\hline
Differential equation method & $0.123285$ & $0.629966$ & $1.00017$ & $1.17647$
& $1.26044$ & $1.29983$\\\hline
Power series & $5.3\times10^{128}$ & $0.629961$ & $1$ & $1.17608$ & $1.25993$
& $1.29928$\\\hline
\end{tabular}
\ \ \
\]

\begin{center}
{\small Table 2: Value of $\sin_{p}\left(  \pi_{p}/2\right)  =\sqrt[p]{p-1}$
obtained independently using each method.}\vspace{0.2cm}
\end{center}

\noindent Notice that the inverse power method appears to be more accurate
when computing $\sin_{p}$ at values close to $\pi_{p}/2$. Indeed, in order to
obtain a good approximation close to this point, it was necessary to allow for
a greater number of terms in the power series than would be necessary for
points far from $\pi_{p}/2$.%

\[%
\begin{tabular}
[c]{|c|c|c|c|c|c|c|}\hline
$p$ & $1.1$ & $1.5$ & $2.0$ & $2.5$ & $3.0$ & $3.5$\\\hline
Inverse power method & $5$ & $8$ & $9$ & $8$ & $8$ & $8$\\\hline
Power Series & $501$ & $13$ & $8$ & $470$ & $501$ & $501$\\\hline
\end{tabular}
\
\]

\begin{center}
{\small Table 3: Number of iterations.}\vspace{0.2cm}\vspace{0.2cm}
\end{center}

\noindent We see that the number of iterations used by the inverse power
method is remarkably low. Below, we present the graphics of $\sin_{p}$ for the
same values of $p$ computed using the three methods (except for $p=1.1$, since
the power series appears to diverge in this case). Notice that all three
methods agree very well with each other, being virtually indistinguishable.

%

\[
\includegraphics[scale=0.4]
{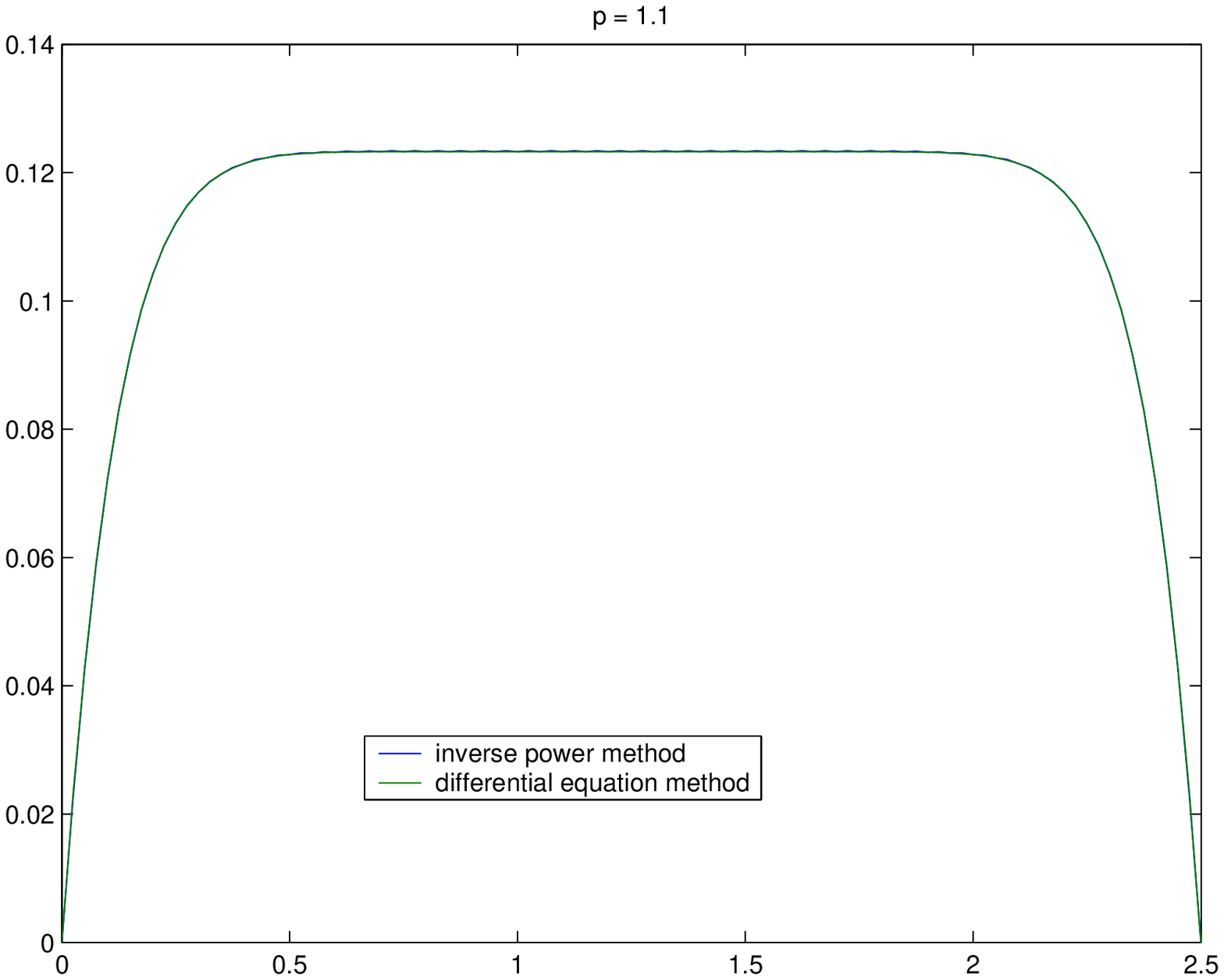}\ \ \ \ \includegraphics[scale=0.4]
{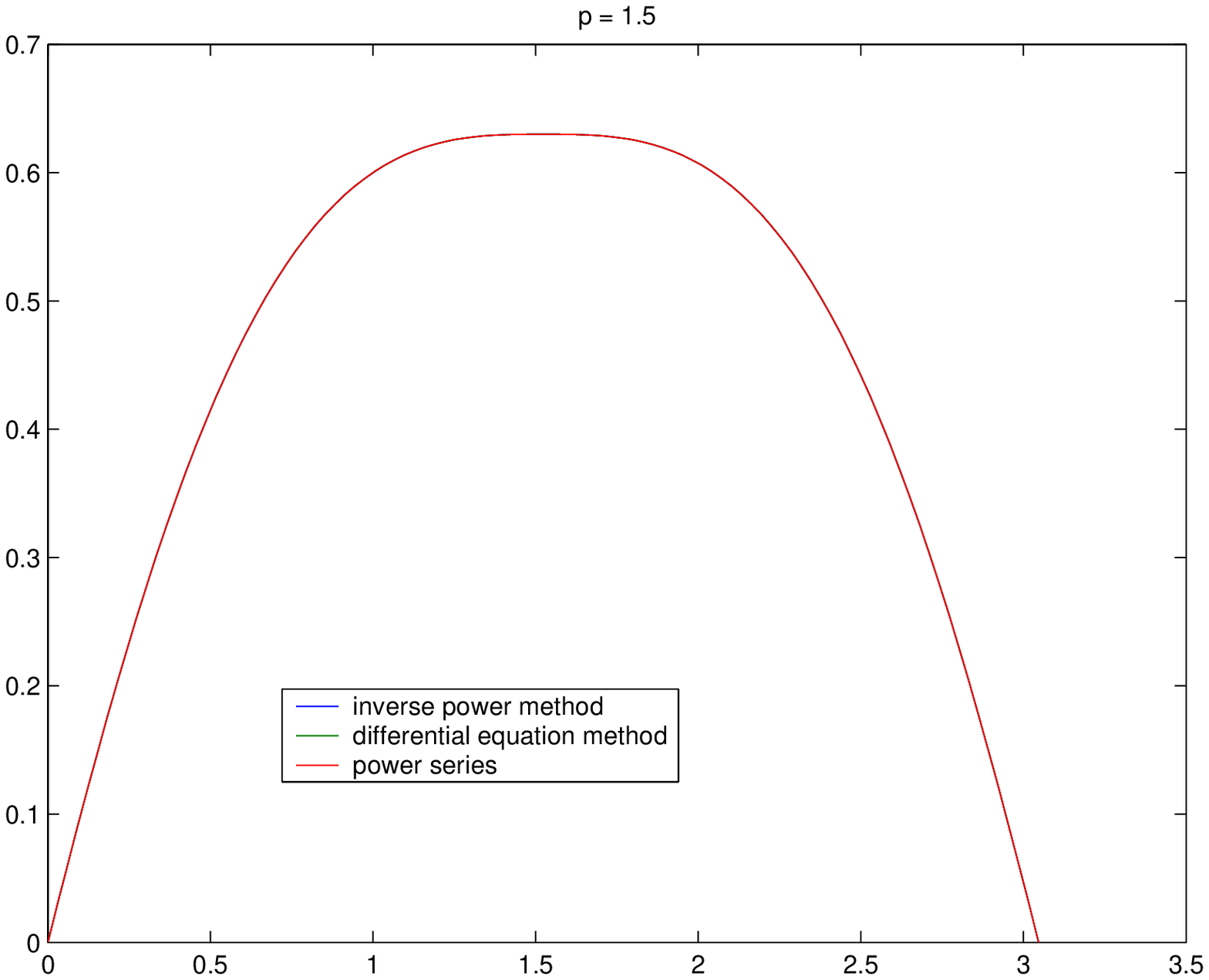}
\]%
\[
\includegraphics[scale=0.4]
{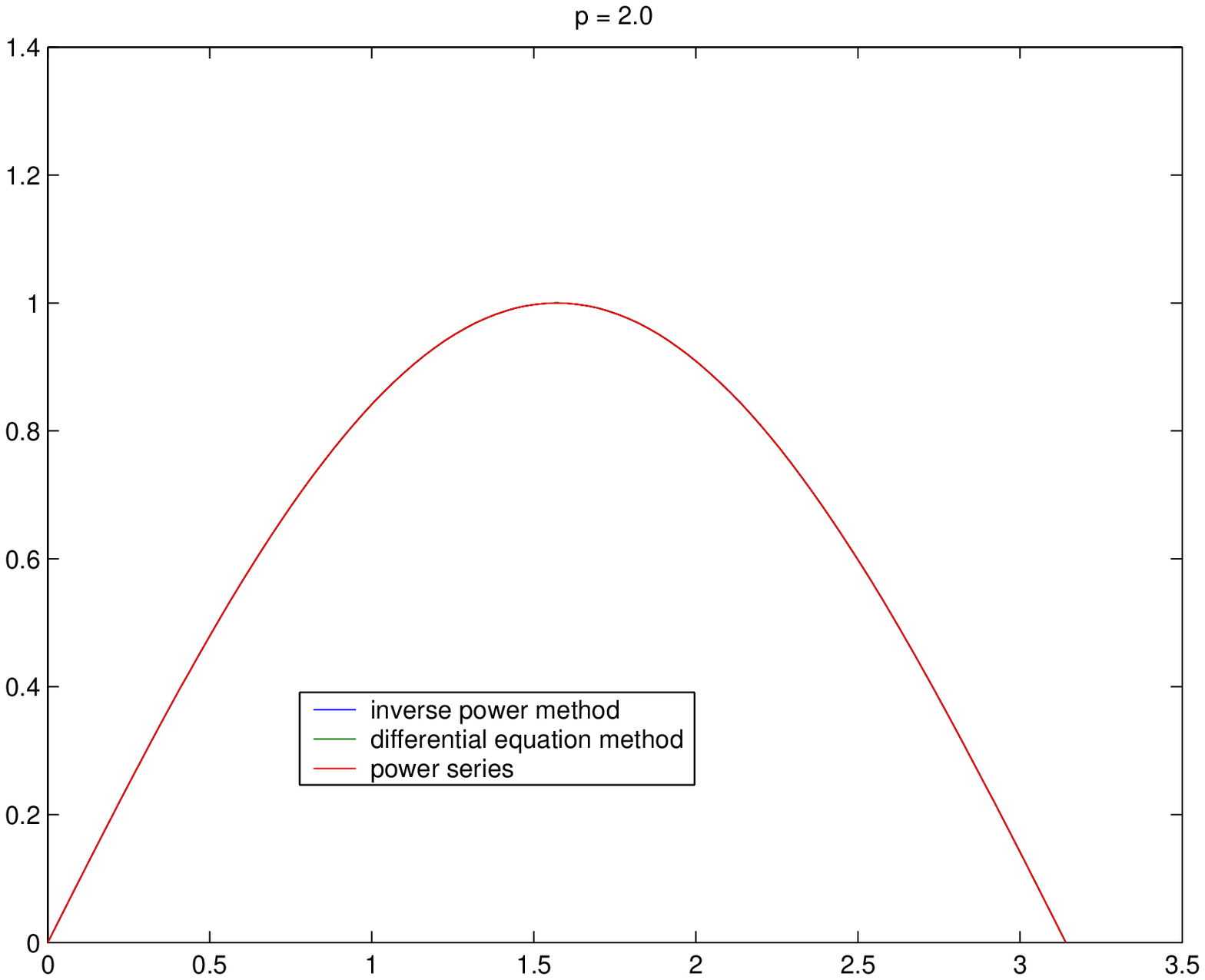}\ \ \ \ \includegraphics[scale=0.4]
{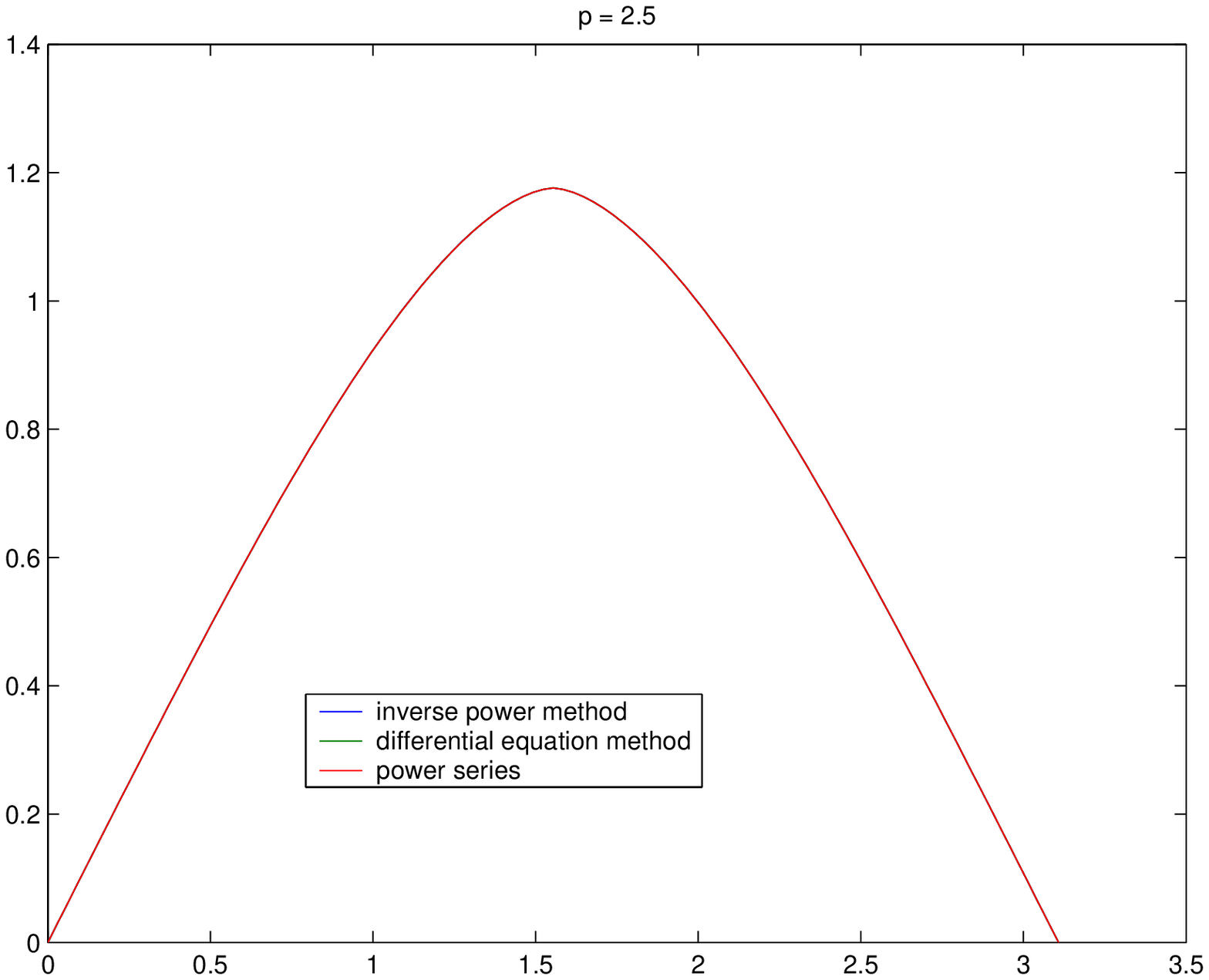}
\]%
\[
\includegraphics[scale=0.4]
{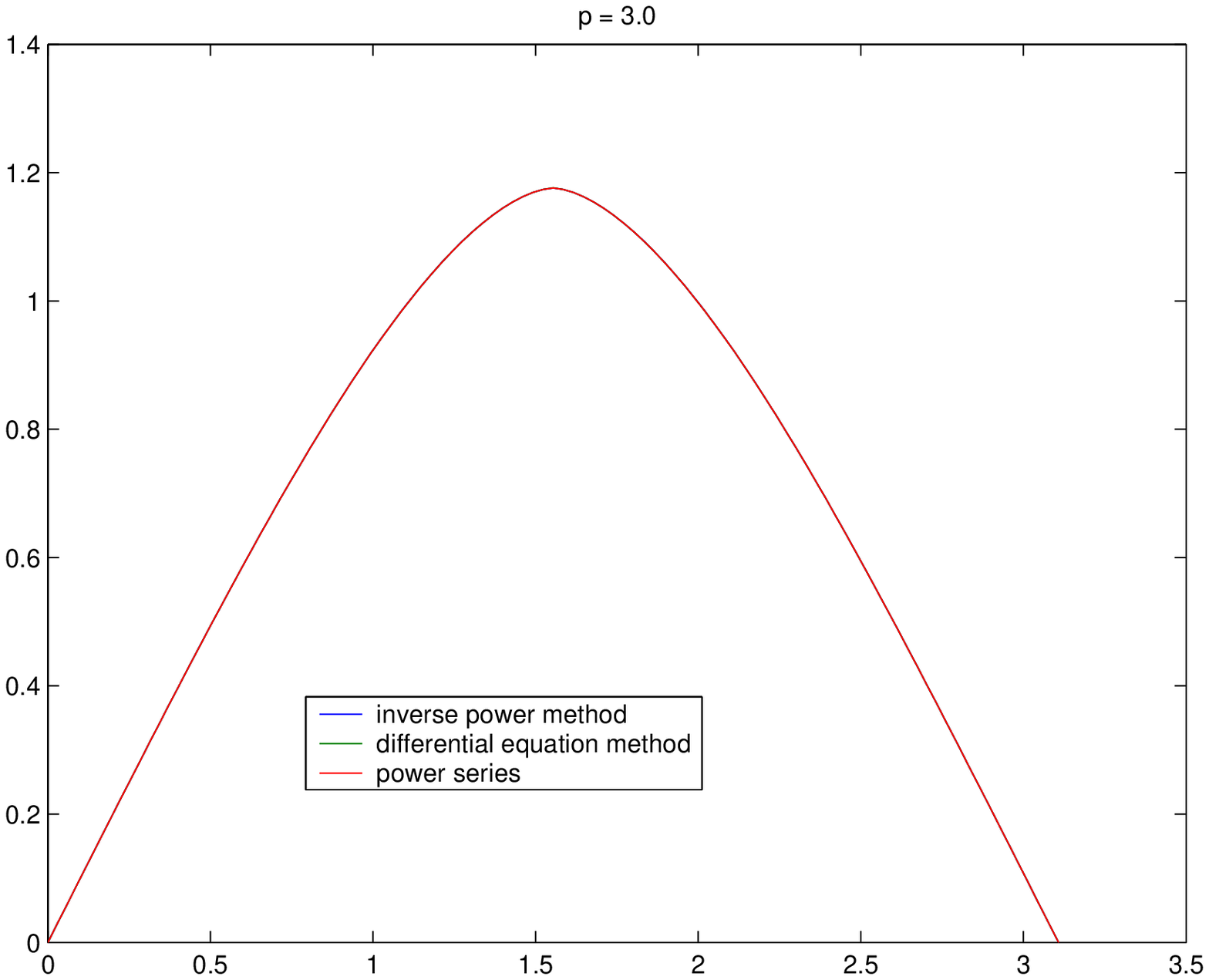}\ \ \ \ \includegraphics[scale=0.4]
{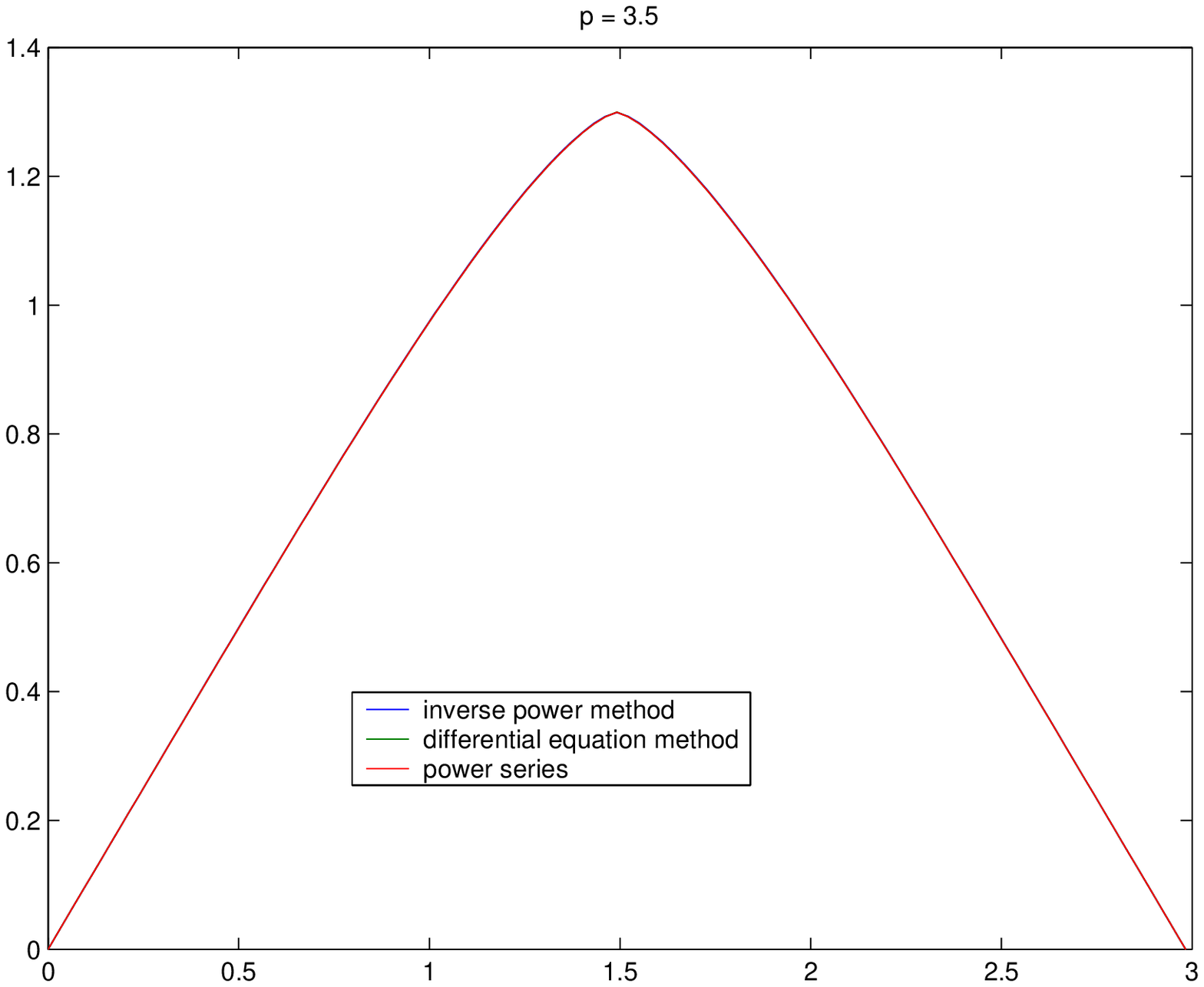}
\]

\section*{Acknowledgments}

The second author would like to thank the support of FAPEMIG and CNPq.

\end{document}